# ENDOGENOUS POST-STRATIFICATION IN SURVEYS: CLASSIFYING WITH A SAMPLE-FITTED MODEL

By F. Jay Breidt[1] and Jean D. Opsomer[2]

*Colorado State University*

Post-stratification is frequently used to improve the precision of survey estimators when categorical auxiliary information is available from sources outside the survey. In natural resource surveys, such information is often obtained from remote sensing data, classified into categories and displayed as pixel-based maps. These maps may be constructed based on classification models fitted to the sample data. Post-stratification of the sample data based on categories derived from the sample data ("endogenous post-stratification") violates the standard post-stratification assumptions that observations are classified without error into post-strata, and post-stratum population counts are known. Properties of the endogenous post-stratification estimator are derived for the case of a sample-fitted generalized linear model, from which the post-strata are constructed by dividing the range of the model predictions into predetermined intervals. Design consistency of the endogenous post-stratification estimator is established under mild conditions. Under a superpopulation model, consistency and asymptotic normality of the endogenous post-stratification estimator are established, showing that it has the same asymptotic variance as the traditional post-stratified estimator with fixed strata. Simulation experiments demonstrate that the practical effect of first fitting a model to the survey data before post-stratifying is small, even for relatively small sample sizes.

**1. Introduction.** Post-stratification (PS) provides a convenient and inexpensive way to improve the precision of estimators in a survey, and is very

Received October 2005; revised March 2007.
[1]Supported in part by USDA Forest Service Rocky Mountain Research Station RJVA 02-JV-11222007-004.
[2]Supported in part by USDA Forest Service Rocky Mountain Research Station RJVA 01-JV-11222007-307.

widely used. In traditional PS, survey observations are classified without error into two or more categories, called post-strata, where the corresponding population counts in those categories are known from some source outside the survey. In surveys of human populations, post-strata are often demographic subgroups, with population counts available from a census. In natural resource surveys, post-strata may be landcover or -use classifications, with population counts obtained from remotely sensed data.

An important example of a natural resource survey is the Forest Inventory and Analysis (FIA) program conducted by the U.S. Forest Service (see, e.g., Frayer and Furnival [5] for a description). In FIA, data of interest are collected annually during intensive field visits and are used to produce official estimates for a large number of forest attributes. In the Interior West region of the United States, FIA estimates are computed as PS estimators, with the strata defined by homogeneous landuse and groundcover categories (e.g., nonforest, broadleaf forest, etc.). Population totals and sample point classifications for those categories are obtained from maps, which are maintained in a geographic information system (GIS). These maps are derived from satellite imagery and other ancillary data layers.

Satellite imagery from the Landsat Enhanced Thematic Mapper Plus (ETM+) as well as from the Moderate Resolution Imaging Spectroradiometer (MODIS) is an important source of remotely sensed data for mapping vegetation over large geographic extents. These data consist of collections of pixel-based maps of physical measurements, such as reflectance values at different wavelengths, which cannot immediately be converted into useable classifications. Instead, categories are obtained by first "training" a classification algorithm on existing satellite imagery and other ancillary digital information, and then predicting the categories of all pixels in the region using that algorithm. Because of the multidimensional and often highly nonlinear nature of the relationships among the variables, the classification algorithms in use today can be quite complex. Examples of such algorithms for forest resources are neural nets and expert systems (Moisen and Frescino [7]). The end result of the classification is a digital (raster) map showing the geographical distribution of the classes over a region of interest. This map is often an important "deliverable" for the organization producing it, and is used by scientists and land managers for a variety of purposes.

Because of the large sample size, detailed nature and high quality of the FIA data, it is attractive to use FIA data to train classification algorithms to produce landcover maps. There are numerous local as well as nationwide mapping efforts that use FIA data for this purpose. Some examples of national efforts include development of the National Landcover Data (http://www.epa.gov/mrlc/nlcd.html), Landfire (http://www.landfire.gov/), and FIA's forest type mapping (Ruefenacht,

Moisen and Blackard [10]). Questions have been raised about the appropriate use of these maps in FIA's PS estimation process (Scott et al. [12]), because the post-strata are delimited with error (since they are based on a model fit) and depend on the sample observations themselves. This violates two fundamental assumptions of traditional PS: the exact post-stratum counts for the population are unknown, and the classification of the sample observations into the post-strata is imperfect. Therefore, it is not clear whether the resulting estimator continues to be consistent and whether the traditional variance estimator remains valid.

We explore the statistical properties of survey estimators that are post-stratified based on a model fitted to the sample observations. To emphasize the relationship between the survey data and the stratification, we will refer to such an estimator as an *endogenous post-stratification estimator*, or EPSE for short. The EPSE is useful in practice whenever population information to construct traditional post-strata is not available, but predictions from a sample-fitted classification model can be generated for the entire population. We restrict our attention to classification schemes based on parametrically specified *generalized linear models* (McCullagh and Nelder [6]). Some of our results will be further restricted to the case of equal-probability sampling, as is used in much of FIA.

An alternative to the EPSE approach is to construct a regression estimator, using the available auxiliary variables as regressors. This can be done using linear models as in the generalized regression estimation (GREG) approach (Cassel, Särndal and Wretman [3]), nonlinear models (Wu and Sitter [14]) or nonparametric models (Breidt and Opsomer [2] and Breidt, Claeskens and Opsomer [1]). Since these models use the auxiliary variables directly, instead of relying only on a classification based on these variables, a properly constructed regression estimator might be more efficient than the EPSE and would have known design properties. However, there are a number of reasons why the EPSE could still be preferable in practice.

First, suitable classification algorithms have already been developed (involving extensive variable selection, model validation and calibration) and maps with well-defined categories are being produced. These maps synthesize information from many layers of geospatial data, so it is operationally efficient to use the generated categories in other estimation problems, rather than building new regression models. Further, categories in the classification can often be readily interpreted (e.g., forest/nonforest), whereas the remote sensing variables (e.g., reflectance at a specific wavelength) are not, so that it is easier to explain the estimation procedure and the resulting fits to diverse end users.

Second, both maps and survey estimates are typically generated by the same organization, so it is clearly desirable to ensure that the survey estimates are calibrated to the map "control" totals. This is automatically achieved under the EPSE approach, but not with a regression estimator.

Third, PS weights based on a modest number of classes may tend to be more stable than the weights obtained from regression estimation, especially in cases where many potentially correlated variables are used in the regression model. In particular, PS weights are guaranteed to be nonnegative, while regression weights are not. Negative weights are an especially serious consideration if survey data are to be used in model fitting, with many statistical programs unable to properly operate in the presence of negative weights.

Finally, the EPSE estimator is robust in the sense that it can compete with the regression estimator when the regression model is correctly specified, and can dominate the regression estimator when the regression model is misspecified.

The EPSE is defined in Section 2 and its properties are described in Section 3, first under a general probability sampling design and then under a superpopulation model. Section 4 describes simulation experiments performed to assess the practical consequences of endogenous PS in a design-based context, and closes with a brief discussion. Proofs are provided in the Appendix.

## 2. Notation and definitions.

2.1. *Post-stratification.* Consider a finite population $U_N = \{1, \ldots, i, \ldots, N\}$. For each $i \in U_N$, an auxiliary vector $\mathbf{x}_i$ is observed. A probability sample $s$ of size $n$ is drawn from $U_N$ according to a sampling design $p_N(\cdot)$, where $p_N(s)$ is the probability of drawing the sample $s$. Assume $\pi_{iN} = \Pr\{i \in s\} = \sum_{s:i \in s} p_N(s) > 0$ for all $i \in U_N$, and define $\pi_{ijN} = \Pr\{i, j \in s\} = \sum_{s:i,j \in s} p_N(s)$ for all $i, j \in U_N$. For compactness of notation we will suppress the subscript $N$ and write $\pi_i$, $\pi_{ij}$ in what follows. Various study variables, generically denoted $y_i$, are observed for $i \in s$.

We now introduce some nonstandard notation for PS that will be useful in our later discussion of endogenous PS. Using the $\{\mathbf{x}_i\}_{i \in U_N}$ and a known vector $\boldsymbol{\lambda}$, a scalar index $\{m(\boldsymbol{\lambda}'\mathbf{x}_i)\}_{i \in U_N}$ is constructed and used to partition $U_N$ into $H$ strata according to predetermined stratum boundaries $-\infty \leq \tau_0 < \tau_1 < \cdots < \tau_{H-1} < \tau_H \leq \infty$. Choice of these boundaries is discussed in Section 3.3 below.

For exponents $\ell = 0, 1, 2$ and stratum indices $h = 1, \ldots, H$, we define

$$A_{Nh\ell}(\boldsymbol{\lambda}) = \frac{1}{N} \sum_{i \in U_N} y_i^\ell I_{\{\tau_{h-1} < m(\boldsymbol{\lambda}'\mathbf{x}_i) \leq \tau_h\}} \tag{1}$$

and

$$A^*_{Nh\ell}(\boldsymbol{\lambda}) = \frac{1}{N} \sum_{i \in U_N} y_i^\ell \frac{I_{\{i \in s\}}}{\pi_i} I_{\{\tau_{h-1} < m(\boldsymbol{\lambda}'\mathbf{x}_i) \leq \tau_h\}}, \tag{2}$$

where $I_{\{C\}} = 1$ if the event $C$ occurs, and zero otherwise. In this notation, stratum $h$ has population stratum proportion $A_{Nh0}(\boldsymbol{\lambda})$, design-weighted sample post-stratum proportion $A^*_{Nh0}(\boldsymbol{\lambda})$, and design-weighted sample post-stratum $y$-mean $A^*_{Nh1}(\boldsymbol{\lambda})/A^*_{Nh0}(\boldsymbol{\lambda})$. The traditional design-weighted PS estimator (PSE) for the population mean $\bar{y}_N = N^{-1} \sum_{i \in U_N} y_i$ is then

$$\hat{\mu}^*_y(\boldsymbol{\lambda}) = \sum_{h=1}^H A_{Nh0}(\boldsymbol{\lambda}) \frac{A^*_{Nh1}(\boldsymbol{\lambda})}{A^*_{Nh0}(\boldsymbol{\lambda})}$$

(3)

$$= \sum_{i \in s} \left\{ \sum_{h=1}^H A_{Nh0}(\boldsymbol{\lambda}) \frac{N^{-1} \pi_i^{-1} I_{\{\tau_{h-1} < m(\boldsymbol{\lambda}' \mathbf{x}_i) \leq \tau_h\}}}{A^*_{Nh0}(\boldsymbol{\lambda})} \right\} y_i = \sum_{i \in s} w^*_{is}(\boldsymbol{\lambda}) y_i,$$

where the sample-dependent weights $\{w^*_{is}(\boldsymbol{\lambda})\}_{i \in s}$ do not depend on $\{y_i\}$, and so can be used for any study variable.

For the important special case of equal-probability designs, in which $\pi_i = nN^{-1}$, we write

$$(4) \qquad A_{nh\ell}(\boldsymbol{\lambda}) = \frac{1}{n} \sum_{i \in s} y_i^\ell I_{\{\tau_{h-1} < m(\boldsymbol{\lambda}' \mathbf{x}_i) \leq \tau_h\}}.$$

In this case, the equal-probability PSE for the population mean $\bar{y}_N$ is

$$(5) \qquad \hat{\mu}_y(\boldsymbol{\lambda}) = \sum_{h=1}^H A_{Nh0}(\boldsymbol{\lambda}) \frac{A_{nh1}(\boldsymbol{\lambda})}{A_{nh0}(\boldsymbol{\lambda})} = \sum_{i \in s} w_{is}(\boldsymbol{\lambda}) y_i,$$

where the weights $\{w_{is}(\boldsymbol{\lambda})\}_{i \in s}$ are obtained by substituting $nN^{-1}$ for $\pi_i$ in (3).

2.2. *Classification based on a generalized linear model.* The notation introduced above does not indicate how the function $m(\cdot)$ might be constructed, nor how values for the parameter vector $\boldsymbol{\lambda}$ should be determined. One possibility is to suppose that a particular study variable, $z_i$, follows a generalized linear model

$$(6) \qquad \mathrm{E}(z_i | \mathbf{x}_i) = m(\boldsymbol{\lambda}' \mathbf{x}_i), \qquad \mathrm{Var}(z_i | \mathbf{x}_i) = v(\mathbf{x}_i),$$

where the expectations are with respect to the model. [For concreteness, think of $z_i$ as a forest/nonforest indicator, with logistic mean function $m(\boldsymbol{\lambda}' \mathbf{x}_i) = \exp(\boldsymbol{\lambda}' \mathbf{x}_i)/\{1 + \exp(\boldsymbol{\lambda}' \mathbf{x}_i)\}$, and $\mathbf{x}_i$ derived from satellite imagery.] We will refer to $z_i$ as the PS variable.

If $\boldsymbol{\lambda}$ were known, we could use $m(\boldsymbol{\lambda}' \mathbf{x}_i)$ as an index to form the PSE in (3) for any study variable $y_i$ (even though the PS is based on a single PS variable $z_i$, the resulting weights can be applied to any response variable $y_i$). If model (6) is true, then $m(\boldsymbol{\lambda}' \mathbf{x}_i)$ is a good predictor for $z_i$. Hence, the estimator (3) applied to the study variable $z_i$ will be more efficient than the Horvitz–Thompson estimator, $\bar{z}_\pi = N^{-1} \sum_{i \in s} \pi_i^{-1} z_i$, which ignores the auxiliary variables $\mathbf{x}_i$. For other study variables $y_i$, the efficiency of (3) relative to $\bar{y}_\pi = N^{-1} \sum_{i \in s} \pi_i^{-1} y_i$ will depend on the relationship between the PS variable $z_i$ and the $y_i$.

2.3. *Endogenous post-stratification.* In endogenous PS, the vector $\boldsymbol{\lambda}$ is unknown, so that estimator (3) is infeasible. Instead, $\boldsymbol{\lambda}$ is estimated from the sample $\{(\mathbf{x}'_i, z_i) : i \in s\}$ by $\hat{\boldsymbol{\lambda}}$ using, for instance, maximum likelihood estimation and for any $i \in U_N$, $z_i$ is predicted by $\hat{z}_i = m(\hat{\boldsymbol{\lambda}}' \mathbf{x}_i)$.

The endogenous post-stratification estimator (EPSE) for the population mean $\bar{y}_N$ is then defined as

$$(7) \qquad \hat{\mu}_y^*(\hat{\boldsymbol{\lambda}}) = \sum_{h=1}^{H} A_{Nh0}(\hat{\boldsymbol{\lambda}}) \frac{A_{Nh1}^*(\hat{\boldsymbol{\lambda}})}{A_{Nh0}^*(\hat{\boldsymbol{\lambda}})} = \sum_{i \in s} w_{is}^*(\hat{\boldsymbol{\lambda}}) y_i.$$

As with the PSE weights, the EPSE weights $\{w_{is}^*(\hat{\boldsymbol{\lambda}})\}_{i \in s}$ as defined in (7) can be applied to any study variable $y$. In the special case of equal-probability designs, we write

$$(8) \qquad \hat{\mu}_y(\hat{\boldsymbol{\lambda}}) = \sum_{h=1}^{H} A_{Nh0}(\hat{\boldsymbol{\lambda}}) \frac{A_{nh1}(\hat{\boldsymbol{\lambda}})}{A_{nh0}(\hat{\boldsymbol{\lambda}})} = \sum_{i \in s} w_{is}(\hat{\boldsymbol{\lambda}}) y_i.$$

Intuitively, it is reasonable to expect that if $\hat{\boldsymbol{\lambda}}$ is a "good" estimator for $\boldsymbol{\lambda}$, then the estimator (7) will behave like the estimator (3), at least asymptotically. We show this equivalence in the sense of design consistency under mild design assumptions in the next section. Such results do not readily yield rates of convergence, because $\hat{\mu}_y^*(\boldsymbol{\lambda})$ is not a differentiably smooth function of $\boldsymbol{\lambda}$, so that traditional Taylor series approaches for the analysis of nonlinear survey estimators (e.g., Särndal et al. [11], Chapter 5) cannot be applied. We therefore restrict our attention to the equal-probability case and study the model-based properties of $\hat{\mu}_y(\hat{\boldsymbol{\lambda}})$, by exploiting the fact that the model expectations of the quantities in (1) and (4) are smooth functions, even though the quantities themselves are not. In particular, we establish a central limit theorem and a consistent variance estimator for $\hat{\mu}_y(\hat{\boldsymbol{\lambda}})$ under an assumed superpopulation model. Section 4 provides simulation evidence that these good model properties also carry over into good design properties.

## 3. Main results.

3.1. *Design assumptions and design consistency.* We assume the general probability sampling design described in Section 2.1 and consider an asymptotic framework in which $N \to \infty$ while the number of strata, $H$, and their boundaries, $\{\tau_h\}$, remain fixed. Assume:

- D1. The covariates $\{\mathbf{x}_i\}$ satisfy $\|\mathbf{x}_i\| \leq M < \infty$. For $\boldsymbol{\lambda} \neq \mathbf{0}$, the empirical distribution function $G_{N\boldsymbol{\lambda}}(z) = N^{-1}\sum_{i\in U_N} I_{\{\mathbf{x}_i'\boldsymbol{\lambda} \leq z\}}$ converges uniformly in $z$ to a limit $G_{\boldsymbol{\lambda}}(z)$, $\lim_{N\to\infty}\sup_z |G_{N\boldsymbol{\lambda}}(z) - G_{\boldsymbol{\lambda}}(z)| = 0$, where the limit is almost sure if the covariates are stochastic.
- D2. The link $m(\cdot)$ is a known, strictly monotone function, $\boldsymbol{\lambda} \neq \mathbf{0}$ is an unknown parameter vector, and $m^{-1}(\tau_h)$ $(h = 1, 2, \ldots, H)$ are continuity points of $G_{\boldsymbol{\lambda}}(z)$. Further, $G_{\boldsymbol{\lambda}}(m^{-1}(\tau_h)) - G_{\boldsymbol{\lambda}}(m^{-1}(\tau_{h-1})) > 0$ for $h = 1, 2, \ldots, H$.
- D3. There is a sequence of estimators of $\boldsymbol{\lambda}$, $\{\hat{\boldsymbol{\lambda}}\}$, with the property that for every $\varepsilon > 0$, there exists $\delta_\varepsilon \in (0, \infty)$ such that $\Pr\{\|\hat{\boldsymbol{\lambda}} - \boldsymbol{\lambda}\| > \delta_\varepsilon\} < \varepsilon$ for all $N$, where the probability is with respect to the sampling design and the covariate model.
- D4. For all $N$, $\min_{i\in U_N} \pi_i \geq \pi_N^* > 0$ where $N\pi_N^* \to \infty$, and there exists $\kappa \geq 0$ such that $N^{1/2+\kappa}(\pi_N^*)^2 \to \infty$ and

$$\max_{i\in U_N} \sum_{j\in U_N: j\neq i} \Delta_{ij}^2 = O(N^{-2\kappa})$$

as $N \to \infty$, where $\Delta_{ij} = \pi_{ij} - \pi_i \pi_j$.
- D5. The study variables $\{y_i\}_{i\in U_N}$ satisfy $\limsup_{N\to\infty} N^{-1}\sum_{i\in U_N} y_i^2 < \infty$.

REMARKS.

1. Note that no stochastic model is assumed for the $\{y_i\}$ in this design-based setting. Randomness comes from the probability mechanism that selects $s$, and possibly from the process generating $\mathbf{x}_i$.
2. The uniform convergence in D1 is met by independent and identically distributed sequences (Glivenko–Cantelli lemma), stationary ergodic sequences (Tucker **(year?)**), certain deterministic sequences [e.g., polynomials of the form $\mathbf{x}_i'\boldsymbol{\lambda} = \sum_{j=0}^p \lambda_j (iN^{-1})^j$], and so forth.
3. D2 ensures that the post-strata are nonempty and can be unambiguously determined from the inverse link; D3 asserts that the parameters in the generalized linear model can be estimated consistently.
4. The first part of D4 allows for sparse sampling in the sense that $\min_{i\in U_N} \pi_i \to 0$ is allowed as $N \to \infty$. The second part of D4 allows for nontrivial dependencies in the sampling. Sparser sampling is possible under weaker design dependence. For example, under simple random sampling without replacement $\max_{i\in U_N} \sum_{j\in U_N: j\neq i} \Delta_{ij}^2 = (N-1)^{-1}(n/N)^2(1-n/N)^2 = O(N^{-1})$ so that D4 holds with $\kappa = 1/2$. On the other hand, consider single-stage cluster sampling of $m$ equally sized clusters from $M$ clusters via simple random sampling without replacement. All elements in each selected cluster are observed. Let $c$ denote the cluster size, and assume it is fixed as $cm = n \to \infty$ and $cM = N \to \infty$. Then

$$\max_{i \in U_N} \sum_{j \in U_N : j \neq i} \Delta_{ij}^2$$

$$= (c-1)\left(\frac{m}{M}\left(1 - \frac{m}{M}\right)\right)^2 + (M-1)c\left(-\frac{m}{M(M-1)}\left(1 - \frac{m}{M}\right)\right)^2$$

$$\leq \left\{(c-1) + \frac{c}{M-1}\right\} = O(1),$$

so that D4 holds with $\kappa = 0$. Note that the corresponding design-covariance assumptions A5 in Robinson and Särndal [9] and A6 in Breidt and Opsomer [2] are not met in general for this design. The Horvitz–Thompson estimator is mean square consistent under D4 and D5, a result of independent interest that is established as a lemma in the Appendix.

RESULT 1. *Assume* D1–D5. *Then the unequal-probability EPSE in* (7) *is design consistent in the sense that for all $\varepsilon > 0$,*

$$\Pr\{|\hat{\mu}_y^*(\hat{\boldsymbol{\lambda}}) - \bar{y}_N| > \varepsilon\} \to 0 \qquad \text{as } N \to \infty.$$

The proof is deferred to the Appendix.

3.2. *Superpopulation model assumptions.* To study further the properties of EPSE, we restrict attention to equal-probability designs and introduce a superpopulation model, which specifies the joint distribution of the random vector $(\mathbf{x}_i', y_i)$, while the randomness of $s$ is not explicitly considered. In what follows, all expectation, probability, and order in probability statements are with respect to this superpopulation model. We continue to consider an asymptotic framework in which $n, N \to \infty$ while the number of strata, $H$, and their boundaries, $\{\tau_h\}$, remain fixed. Our proofs rely on the approach of Randles [8]. Formally, we assume the following:

- M1. The covariates $\{\mathbf{x}_i\}$ are independent and identically distributed (i.i.d.) random $p$-vectors with nondegenerate continuous joint probability density function $f$ and compact support.
- M2. The link $m(\cdot)$ is a known, strictly monotone function on its domain, $\boldsymbol{\lambda} \neq \mathbf{0}$ is an unknown parameter vector and $v(\cdot)$ in (6) is a bounded, positive function.
- M3. There is a sequence of estimators of $\boldsymbol{\lambda}$, $\{\hat{\boldsymbol{\lambda}}\}$, such that $\hat{\boldsymbol{\lambda}} - \boldsymbol{\lambda} = O_p(n^{-1/2})$.
- M4. The study variables $y_i | \mathbf{x}_i$ are conditionally independent random variables with $\mathrm{E}(y_i^4 | \mathbf{x}_i) \leq K_1 < \infty$. Also,

$$\alpha_{h\ell}(\boldsymbol{\lambda}) = \mathrm{E}(y_i^\ell I_{\{\tau_{h-1} < m(\boldsymbol{\lambda}'\mathbf{x}_i) \leq \tau_h\}}) \tag{9}$$

is continuous in $\boldsymbol{\lambda}$ for $\ell = 0, 1, 2$, and $\alpha_{h0}(\boldsymbol{\lambda}) > 0$ for $h = 1, \ldots, H$. In particular, the variables $z_i | \mathbf{x}_i$ are conditionally independent random variables with

$$\mathrm{E}(z_i | \mathbf{x}_i) = m(\boldsymbol{\lambda}' \mathbf{x}_i), \qquad \mathrm{Var}(z_i | \mathbf{x}_i) = v(\mathbf{x}_i), \qquad \mathrm{E}(z_i^4 | \mathbf{x}_i) \leq K_1 < \infty.$$

- M5. The sample $s$ is selected according to an equal-probability design of fixed size $n$, with $\pi_i = nN^{-1} \to \pi \in [0, 1]$ as $n, N \to \infty$.

While the conditional independence in M4 rules out certain clustered designs, it seems quite plausible in large-scale natural resource surveys, where it is often the case that sampling locations are widely dispersed and, after correcting for covariates, no spatial dependence remains. (We investigate the effect of residual spatial dependence via simulation in Section 4.) The equal-probability design assumed in M5 is also somewhat limiting, but it does cover the systematic designs used by the U.S. Forest Service in FIA. Further, our results extend trivially to the case of a fixed number of design strata (determined prior to sampling, unlike post-strata) with a large equal-probability sample within each stratum, and possibly unequal probabilities across strata.

3.3. *Central limit theorem.* The proof of consistency and asymptotic normality for the EPSE in (8) with respect to the superpopulation model is deferred to the Appendix.

RESULT 2. *Under assumptions* M1–M5,

$$\left\{ \frac{1}{n} \left(1 - \frac{n}{N}\right) \right\}^{-1/2} (\hat{\mu}_y(\hat{\boldsymbol{\lambda}}) - \bar{y}_N) \xrightarrow{\mathcal{L}} \mathcal{N}(0, V_{y\boldsymbol{\lambda}}),$$

*where*

$$V_{y\boldsymbol{\lambda}} = \sum_{h=1}^{H} \Pr\{\tau_{h-1} < m(\boldsymbol{\lambda}' \mathbf{x}_i) \leq \tau_h\} \mathrm{Var}(y_i | \tau_{h-1} < m(\boldsymbol{\lambda}' \mathbf{x}_i) \leq \tau_h).$$

REMARKS.

1. When $H = 1$, the asymptotic model variance of $\hat{\mu}_y(\hat{\boldsymbol{\lambda}}) - \bar{y}_N$ from Result 2 is

(10) $$\frac{1}{n}\left(1 - \frac{n}{N}\right) \mathrm{Var}(y_i).$$

This is the model variance of $\bar{y}_\pi - \bar{y}_N$ under any equal-probability design, or the model-averaged design variance of the ordinary sample mean under simple random sampling without replacement.

2. For general $H$, the asymptotic model variance in Result 2 is equal to that of the traditional post-stratified estimator (i.e., in which $\boldsymbol{\lambda}$ is known). This variance has an intuitive form: it sums stratum fraction times within-stratum variance of $y$ over the post-strata. Note that

$$\sum_{h=1}^{H} \alpha_{h0}(\boldsymbol{\lambda})\left(\frac{\alpha_{h1}(\boldsymbol{\lambda})}{\alpha_{h0}(\boldsymbol{\lambda})}\right)^2 \geq \left(\sum_{h=1}^{H} \alpha_{h0}(\boldsymbol{\lambda})\frac{\alpha_{h1}(\boldsymbol{\lambda})}{\alpha_{h0}(\boldsymbol{\lambda})}\right)^2 = (\mathrm{E}[y_i])^2,$$

with equality only if the post-stratum means are all identical: $\alpha_{h1}(\boldsymbol{\lambda})/\alpha_{h0}(\boldsymbol{\lambda}) \equiv \mathrm{E}(y_i)$ for $h=1,\ldots,H$. Thus, we have that

$$\mathrm{Var}(y_i) = \sum_{h=1}^{H} \alpha_{h0}(\boldsymbol{\lambda})\frac{\alpha_{h2}(\boldsymbol{\lambda})}{\alpha_{h0}(\boldsymbol{\lambda})} - (\mathrm{E}[y_i])^2$$

$$\geq \sum_{h=1}^{H} \alpha_{h0}(\boldsymbol{\lambda})\left(\frac{\alpha_{h2}(\boldsymbol{\lambda})}{\alpha_{h0}(\boldsymbol{\lambda})} - \left(\frac{\alpha_{h1}(\boldsymbol{\lambda})}{\alpha_{h0}(\boldsymbol{\lambda})}\right)^2\right)$$

$$= \sum_{h=1}^{H} \Pr\{\tau_{h-1} < m(\boldsymbol{\lambda}'\mathbf{x}_i) \leq \tau_h\}\mathrm{Var}(y_i|\tau_{h-1} < m(\boldsymbol{\lambda}'\mathbf{x}_i) \leq \tau_h),$$

so that unless the post-stratum means of $y$ are all identical, the EPSE will be asymptotically more efficient than the Horvitz–Thompson estimator (the ordinary sample mean in this case).

3. For the PS variable $z_i$, it can be shown that

$$V_{z\boldsymbol{\lambda}} = \mathrm{E}[v(\mathbf{x}_i)] + \sum_{h=1}^{H} \Pr\{\tau_{h-1} < m(\boldsymbol{\lambda}'\mathbf{x}_i) \leq \tau_h\} \tag{11}$$

$$\times \mathrm{Var}(m(\boldsymbol{\lambda}'\mathbf{x}_i)|\tau_{h-1} < m(\boldsymbol{\lambda}'\mathbf{x}_i) \leq \tau_h)$$

with $v(\mathbf{x}_i) = \mathrm{Var}(z_i|\mathbf{x}_i)$. A lower bound for the asymptotic model variance of $\hat{\mu}_z(\hat{\boldsymbol{\lambda}}) - \bar{z}_N$ is given by

$$\frac{1}{n}\left(1 - \frac{n}{N}\right)\mathrm{E}[v(\mathbf{x}_i)],$$

which is also the asymptotic model variance of the error of the nonlinear regression estimator

$$\hat{\eta}_z(\hat{\boldsymbol{\lambda}}) = N^{-1}\sum_{i \in U_N} m(\hat{\boldsymbol{\lambda}}'\mathbf{x}_i) + n^{-1}\sum_{i \in s}(z_i - m(\hat{\boldsymbol{\lambda}}'\mathbf{x}_i)). \tag{12}$$

Hence, the quantity

$$\frac{1}{n}\left(1 - \frac{n}{N}\right)(V_{z\boldsymbol{\lambda}} - \mathrm{E}[v(\mathbf{x}_i)])$$

measures the asymptotic loss in efficiency of the EPSE relative to the regression estimator $\hat{\eta}_z(\hat{\boldsymbol{\lambda}})$. The EPSE will be as asymptotically efficient as the regression estimator (12), that is, $V_{z\boldsymbol{\lambda}} = \mathrm{E}[v(\mathbf{x}_i)]$, if the $m(\boldsymbol{\lambda}'\mathbf{x}_i)$ are constant within each stratum. If this is not the case, the EPSE will fail to match the asymptotic efficiency of the regression estimator. It should be noted that although the EPSE for the study variable $z$ is therefore likely to be dominated by the regression estimator (12), this is not necessarily the case if a different regression estimator is used. For other study variables $y$, the EPSE may be better than a regression estimator, depending on the relationship between $y$ and $\mathbf{x}$ in the population. We explore this further in the simulation study in the next section.

4. In some applications, it might be of interest to select the $\{\tau_h\}$ defining the categories to improve the efficiency of the EPSE for a "target" variable $z$. As noted in Section 1, these class boundaries are often determined by the requirements of the classification algorithm and the desired map output to which the EPSE is calibrated, so that little choice might be available when they are applied in the construction of the post-strata, except for possibly collapsing neighboring post-strata in case of small sample sizes. If the operational environment allows for the selection of stratum boundaries, then boundaries might be constructed by applying the cumulative root-density method described in Cochran [4], Section 5A.7, to the $m(\hat{\boldsymbol{\lambda}}'\mathbf{x}_i)$, though this method requires further study.

3.4. *Variance estimation.* We now consider variance estimation for the EPSE. The standard post-stratified design variance estimator under simple random sampling without replacement is

$$
\hat{V}(\hat{\mu}_y(\boldsymbol{\lambda})) = \frac{1}{n}\left(1 - \frac{n}{N}\right) \sum_{h=1}^{H} \frac{(N_h N^{-1})^2}{n_h n^{-1}} s_{yh}^2
$$
(13)
$$
= \frac{1}{n}\left(1 - \frac{n}{N}\right)\left\{\sum_{h=1}^{H} \frac{A_{Nh0}^2(\boldsymbol{\lambda})}{A_{nh0}(\boldsymbol{\lambda})} \frac{A_{nh2}(\boldsymbol{\lambda}) - A_{nh1}^2(\boldsymbol{\lambda})/A_{nh0}(\boldsymbol{\lambda})}{A_{nh0}(\boldsymbol{\lambda}) - n^{-1}}\right\},
$$

where $N_h$ and $n_h$ are population and sample counts within post-stratum $h$, and $s_{yh}^2$ is the sample $y$-variance within post-stratum $h$ (see, e.g., Särndal et al. [11], equation (7.6.5)).

The next result shows that the analogous estimator under endogenous PS consistently estimates the asymptotic model variance of the EPSE. The proof is again deferred to the Appendix.

RESULT 3. *Let*

(14) $\quad \hat{V}(\hat{\mu}_y(\hat{\boldsymbol{\lambda}})) = \frac{1}{n}\left(1 - \frac{n}{N}\right)\left\{\sum_{h=1}^{H} \frac{A_{Nh0}^2(\hat{\boldsymbol{\lambda}})}{A_{nh0}(\hat{\boldsymbol{\lambda}})} \frac{A_{nh2}(\hat{\boldsymbol{\lambda}}) - A_{nh1}^2(\hat{\boldsymbol{\lambda}})/A_{nh0}(\hat{\boldsymbol{\lambda}})}{A_{nh0}(\hat{\boldsymbol{\lambda}}) - n^{-1}}\right\}.$

*Under assumptions* M1–M5, *as* $n, N \to \infty$

$$\hat{V}(\hat{\mu}_y(\hat{\boldsymbol{\lambda}}))^{-1/2}(\hat{\mu}_y(\hat{\boldsymbol{\lambda}}) - \bar{y}_N) \xrightarrow{\mathcal{L}} N(0,1).$$

**4. Simulations.** The two main goals of the simulation are to assess the design efficiency of the EPSE and the design bias of the variance estimator (14); we also look at confidence interval coverage. The simulations are performed in a setting that mimics a real survey, in which characteristics of multiple study variables are estimated using one set of weights. The weights for estimation of a mean are the Horvitz–Thompson estimator (HTE) weights $\{n^{-1}\}_{i \in s}$, the PSE weights $\{w_{is}(\boldsymbol{\lambda})\}_{i \in s}$, the EPSE weights $\{w_{is}(\hat{\boldsymbol{\lambda}})\}_{i \in s}$, or the simple linear regression (REG) weights (e.g., (6.5.12) of Särndal et al. [11]). The HTE does not use auxiliary information; the PSE uses auxiliary information with a known model; and the EPSE and REG use auxiliary information with fitted models.

We consider two different models for the PS variable, $z_i$. First, we look at the situation in which the true model in (6) for the PS variable is continuous and follows a ratio model (see, e.g., Särndal et al. [11], page 226), so that $m(\cdot)$ is the identity function. Second, we consider the case where the PS variable is binary and $m(\cdot)$ is the logistic link.

4.1. *Ratio model post-stratification.* We first describe the simulation setup for the ratio model PS. We assume a population of size $N = 1000$ with eight survey variables of interest. For the PS variable $z_i$, we let $\mathrm{E}(z_i|\mathbf{x}) = m(\boldsymbol{\lambda}'\mathbf{x}) = 1 + 2(x - 0.5)$ and $\mathrm{Var}(z_i|\mathbf{x}) = v(x) = 2\sigma^2 x$, while for the remaining seven variables $(y_i)$, we take their mean functions to be equal to $g_k$,

$$\begin{aligned}
\text{quadratic:} \quad & g_1(x) = 1 + 2(x - 0.5)^2, \\
\text{bump:} \quad & g_2(x) = 1 + 2(x - 0.5) + \exp(-200(x - 0.5)^2), \\
\text{jump:} \quad & g_3(x) = \{1 + 2(x - 0.5)\} I_{\{x \le 0.65\}} + 0.65 I_{\{x > 0.65\}}, \\
\text{exponential:} \quad & g_4(x) = \exp(-8x), \\
\text{cycle 1:} \quad & g_5(x) = 2 + \sin(2\pi x), \\
\text{cycle 2:} \quad & g_6(x) = 2 + \sin(8\pi x), \\
\text{white noise:} \quad & g_7(x) = 8,
\end{aligned}$$

and variance equal to $\sigma^2$, with $x$ uniformly distributed on $(0, 1)$ and the errors for all functions independent and normally distributed. The variance function for the PS variable is chosen so that, averaging over the covariate $x$, we have $\mathrm{E}[v(x)] = \sigma^2$. Thus, the PS variable and the remaining seven study variables all have the same variance, averaged over $x$. We considered two different values of $\sigma$: 0.25 and 0.50.

For each noise level, we fixed the population (i.e., simulated $N$ values for each of the eight variables of interest) and drew 1000 replicate samples of size $n = 50$ and $n = 100$, each via simple random sampling without replacement from this fixed population. We constructed HTE and REG weights using standard methods. We used the model for the PS variable for constructing the PSE weights, using known parameter values, and the EPSE weights, using fitted parameter values. The weights were then applied to the remaining seven study variables. Hence, the PS model will only be correctly specified for the PS variable. The EPSE and PSE were calculated using two strata with boundaries $\tau = (-\infty, 1.0, \infty)$, using four strata with boundaries $\tau = (-\infty, 0.5, 1.0, 1.5, \infty)$, and using six strata with boundaries $\tau = (-\infty, 0.5, 0.75, 1.0, 1.25, 1.5, \infty)$.

Table 1 summarizes the design efficiency results as ratios of the MSE of the HTE, PSE($H$), or REG over the MSE of the EPSE($H$), where $H = 2, 4$ or 6 strata. Overall, the results show that the EPSE behaves as expected: it produces a large improvement in efficiency relative to the HTE for the variable on which the PS is based, as well as for most of the other variables that are correlated with it. When the number of strata increases, the efficiency gains become more pronounced, though EPSE begins to break down due to post-stratum sample sizes of zero or 1 when the number of strata is large and the sample size is small. When the relationship between the variables of interest and the auxiliary variable becomes less strong (i.e., higher noise levels), the efficiency gains of EPSE decrease. EPSE is typically as good as or better than REG for study variables on which the regression model is badly misspecified, but loses out to REG when the true model is linear or nearly so ("bump"). The "white noise" variable shows that, when a variable is not related to the stratification variable, the efficiency is near that of the HTE, but with decreasing efficiency as the number of strata increases (since the strata are entirely unnecessary).

Table 1 also shows that the EPSE is essentially equivalent to the PSE in terms of design efficiency, even for $n = 50$, implying that the effect of basing the PS on a fitted model instead of on exogenous strata is negligible for moderate to large sample sizes.

Next, we consider the variance estimator proposed in (14) by computing percentage relative biases (100% times the variance bias divided by the true design variance) for the PSE variance estimator (13) and the EPSE variance estimator (14). These results (not tabled) show that neither estimator is unbiased, and both tend to show negative bias (147 of the 192 cases in Table 1). The bias of the EPSE variance estimator tracks that of the PSE variance estimator closely for low noise, low number of strata and large sample size, but the tracking deteriorates as noise increases, number of strata increases or sample size decreases.

Finally, we assess the quality of the normal approximation by constructing approximate 95% confidence intervals from the pivotal quantity in Result 3. These confidence intervals, $\hat{\mu}_y(\hat{\boldsymbol{\lambda}}) \pm 1.96\{\hat{V}(\hat{\mu}_y(\hat{\boldsymbol{\lambda}}))\}^{1/2}$, attained empirical

TABLE 1
*Ratio of MSE of Horvitz–Thompson (HTE), post-stratification on H strata [PSE(H)] and linear regression (REG) estimators to MSE of endogenous post-stratification estimator on H strata [EPSE(H)]*

| Response Variable | $\sigma$ | $n$ | EPSE(2) versus | | | EPSE(4) versus | | | EPSE(6) versus | | |
|---|---|---|---|---|---|---|---|---|---|---|---|
| | | | HTE | PSE(2) | REG | HTE | PSE(4) | REG | HTE | PSE(6) | REG |
| PS variable | 0.25 | 50 | 2.76 | 1.06 | 0.42 | 4.57 | 1.02 | 0.70 | 4.95 | 0.99 | 0.73 |
| | 0.25 | 100 | 2.46 | 1.01 | 0.41 | 4.36 | 1.01 | 0.72 | 4.75 | 1.01 | 0.78 |
| | 0.50 | 50 | 1.78 | 1.04 | 0.74 | 2.14 | 1.01 | 0.87 | 2.17 | 1.01 | 0.86 |
| | 0.50 | 100 | 1.65 | 1.04 | 0.74 | 2.01 | 1.02 | 0.90 | 2.00 | 1.01 | 0.90 |
| quad | 0.25 | 50 | 0.94 | 1.00 | 1.01 | 1.09 | 1.01 | 1.18 | 1.04 | 0.98 | 1.10 |
| | 0.25 | 100 | 0.93 | 1.00 | 1.00 | 1.10 | 1.01 | 1.19 | 1.08 | 0.99 | 1.16 |
| | 0.50 | 50 | 0.95 | 1.00 | 1.00 | 0.94 | 1.00 | 1.00 | 0.89 | 0.99 | 0.94 |
| | 0.50 | 100 | 0.93 | 1.00 | 0.99 | 0.96 | 0.99 | 1.02 | 0.93 | 0.98 | 1.00 |
| bump | 0.25 | 50 | 2.17 | 0.98 | 0.70 | 3.13 | 0.97 | 1.01 | 4.29 | 0.99 | 1.35 |
| | 0.25 | 100 | 2.16 | 0.98 | 0.69 | 3.21 | 0.97 | 1.02 | 4.20 | 0.98 | 1.33 |
| | 0.50 | 50 | 1.55 | 0.97 | 0.83 | 1.84 | 0.99 | 0.98 | 2.06 | 1.00 | 1.10 |
| | 0.50 | 100 | 1.57 | 0.98 | 0.84 | 1.86 | 0.99 | 1.01 | 2.09 | 1.00 | 1.13 |
| jump | 0.25 | 50 | 1.09 | 0.99 | 0.99 | 1.54 | 0.99 | 1.40 | 1.79 | 0.97 | 1.60 |
| | 0.25 | 100 | 1.10 | 0.99 | 1.00 | 1.55 | 0.98 | 1.41 | 1.83 | 0.99 | 1.66 |
| | 0.50 | 50 | 1.01 | 1.00 | 0.99 | 1.13 | 1.01 | 1.11 | 1.12 | 0.97 | 1.10 |
| | 0.50 | 100 | 1.00 | 0.99 | 1.00 | 1.14 | 1.00 | 1.14 | 1.17 | 0.98 | 1.17 |
| expon | 0.25 | 50 | 1.12 | 1.01 | 0.88 | 1.26 | 1.01 | 1.00 | 1.19 | 0.99 | 0.94 |
| | 0.25 | 100 | 1.08 | 1.00 | 0.85 | 1.31 | 1.00 | 1.04 | 1.29 | 0.99 | 1.02 |
| | 0.50 | 50 | 1.02 | 1.00 | 0.96 | 1.01 | 1.01 | 0.95 | 0.94 | 0.99 | 0.90 |
| | 0.50 | 100 | 0.98 | 1.01 | 0.94 | 1.03 | 0.99 | 0.98 | 1.01 | 0.98 | 0.96 |

TABLE 1
*(Continued)*

| Response Variable | $\sigma$ | $n$ | EPSE(2) versus | | | EPSE(4) versus | | | EPSE(6) versus | | |
|---|---|---|---|---|---|---|---|---|---|---|---|
| | | | HTE | PSE(2) | REG | HTE | PSE(4) | REG | HTE | PSE(6) | REG |
| cycle 1 | 0.25 | 50 | 3.59 | 1.02 | 1.67 | 3.57 | 1.04 | 1.66 | 4.84 | 1.00 | 2.09 |
| | 0.25 | 100 | 3.77 | 1.00 | 1.65 | 3.73 | 1.00 | 1.64 | 4.82 | 1.01 | 2.11 |
| | 0.50 | 50 | 2.11 | 0.98 | 1.28 | 2.14 | 1.04 | 1.30 | 2.33 | 0.99 | 1.35 |
| | 0.50 | 100 | 2.23 | 1.01 | 1.28 | 2.23 | 1.02 | 1.29 | 2.43 | 1.00 | 1.40 |
| cycle 4 | 0.25 | 50 | 1.05 | 1.00 | 0.98 | 1.00 | 1.00 | 0.93 | 1.40 | 0.95 | 1.28 |
| | 0.25 | 100 | 1.07 | 1.00 | 0.97 | 1.05 | 1.01 | 0.95 | 1.59 | 0.98 | 1.45 |
| | 0.50 | 50 | 1.02 | 1.01 | 0.98 | 0.97 | 1.00 | 0.93 | 1.12 | 0.89 | 1.06 |
| | 0.50 | 100 | 1.06 | 1.01 | 0.98 | 1.04 | 1.02 | 0.96 | 1.26 | 0.89 | 1.16 |
| white noise | 0.25 | 50 | 0.96 | 1.00 | 1.00 | 0.91 | 1.00 | 0.94 | 0.86 | 0.98 | 0.89 |
| | 0.25 | 100 | 0.94 | 1.00 | 0.99 | 0.92 | 1.00 | 0.98 | 0.90 | 0.99 | 0.96 |
| | 0.50 | 50 | 0.96 | 1.00 | 1.00 | 0.90 | 1.00 | 0.94 | 0.85 | 0.98 | 0.89 |
| | 0.50 | 100 | 0.94 | 1.00 | 0.99 | 0.92 | 0.99 | 0.97 | 0.90 | 0.98 | 0.95 |

Numbers greater than 1 favor EPSE. Based on ratio model post-stratification in 1000 replications of simple random sampling from a fixed population of size $N = 1000$. Replications in which at least one stratum had fewer than two samples are omitted from the summary: 55 reps for six strata at $n = 50$, $\sigma = 0.25$; 58 reps for six strata at $n = 50$, $\sigma = 0.5$; and 3 reps for four strata at $n = 50$, $\sigma = 0.5$.

coverages (not tabled) ranging from 92.1% to 95.8% for the 96 combinations of noise level, sample size, number of strata and study variable. These empirical coverages track closely the empirical coverages of confidence intervals constructed from the PSE.

We repeated the experiments with $\sigma = 0.25$, $n = 50$, and $H = 2, 4$ or 6 strata for the case with $x_i = i/(N+1)$ $(i = 1, 2, \ldots, N)$ and with the residuals for every variable autocorrelated: $\text{corr}(z_i, z_j) = \text{corr}(y_i, y_j) = 0.99^{|i-j|}$. This setting clearly violates the conditional independence assumption of M4. The results (not tabled) indicate that the EPSE remains essentially unbiased and its confidence intervals continue to have close to nominal coverage (92.9%–95.6%). Efficiency compared to HTE is even greater than in the conditionally independent case, because positive autocorrelation is trend-like behavior that is captured to some extent by the post-strata. The variance estimator tends to have less negative bias or more positive bias than in the conditional independence case. Overall, these limited simulations suggest that EPSE maintains its good behavior outside of the limited setting described in the technical assumptions of Section 3.

4.2. *Logistic model post-stratification.* Since the theory of this paper covers generalized linear models, the above simulation experiments were repeated after replacing the ratio model for the PS variable by a logistic model, but keeping all other models the same. In this case, the PS variable $z_i$ is now a Bernoulli variable with expectation $m(x) = \exp(\lambda_0 + \lambda_1 x)/(1 + \exp(\lambda_0 + \lambda_1 x))$. The values for $(\lambda_0, \lambda_1)$ were chosen as $(-10, 20)$ for the "low noise" case ($\sigma = 0.25$ for the remaining variables) and as $(-3, 6)$ for the "high noise" case ($\sigma = 0.50$ for the other variables). These levels will be denoted as the "steep" and the "flat" model. Two, four and six equal-size strata partitioning $[0, 1]$ are considered for the PSE and EPSE. All estimators remain as in Section 4.1.

Table 2 displays the relative efficiency for the logistic model simulations using $n = 200$ (the logistic fits were problematic at smaller sample sizes). The findings are very similar to those discussed for the ratio model. The EPSE continues to improve substantially over the HTE for most variables, while not deviating substantially from the PSE with known stratum classifications. Further, the EPSE continues to be competitive with the REG estimator. Efficiency tends to increase from two to four strata, but level off from four to six strata.

Approximate 95% confidence intervals computed from the pivotal quantity in Result 3 attained empirical coverages (not tabled) ranging from 93.7% to 96.3% for the 48 combinations of model, number of strata and study variable in Table 2. These empirical coverages track closely the empirical coverages of confidence intervals constructed from the PSE, and are quite close to nominal in spite of finite-sample bias in the variance estimator.

ENDOGENOUS POST-STRATIFICATION 17

TABLE 2
*Ratio of MSE of Horvitz–Thompson (HTE), post-stratification on H strata [PSE(H)] and linear regression (REG) estimators to MSE of endogenous post-stratification estimator on H strata [EPSE(H)]*

| Population | | $\sigma$ | EPSE(2) versus | | | EPSE(4) versus | | | EPSE(6) versus | | |
|---|---|---|---|---|---|---|---|---|---|---|---|
| | | | HTE | PSE(2) | REG | HTE | PSE(4) | REG | HTE | PSE(6) | REG |
| PS | steep | 0.25 | 3.62 | 1.01 | 1.15 | 4.33 | 0.93 | 1.37 | 4.59 | 0.98 | 1.44 |
| variable | flat | 0.50 | 1.42 | 1.02 | 0.88 | 1.52 | 1.00 | 0.95 | 1.51 | 0.98 | 0.94 |
| quad | steep | 0.25 | 1.07 | 1.00 | 1.00 | 1.13 | 1.01 | 1.06 | 1.14 | 1.00 | 1.06 |
| | flat | 0.50 | 1.08 | 1.00 | 1.00 | 1.12 | 1.01 | 1.04 | 1.10 | 0.99 | 1.02 |
| bump | steep | 0.25 | 2.32 | 0.94 | 0.61 | 3.46 | 0.95 | 0.91 | 3.85 | 0.94 | 1.03 |
| | flat | 0.50 | 1.73 | 0.96 | 0.77 | 2.22 | 1.00 | 0.99 | 2.49 | 0.99 | 1.11 |
| jump | steep | 0.25 | 1.11 | 0.99 | 0.97 | 1.21 | 0.99 | 1.05 | 1.32 | 1.01 | 1.15 |
| | flat | 0.50 | 1.05 | 0.99 | 0.98 | 1.25 | 0.97 | 1.16 | 1.29 | 1.04 | 1.20 |
| expon | steep | 0.25 | 1.25 | 0.99 | 0.87 | 1.27 | 1.00 | 0.89 | 1.29 | 1.00 | 0.89 |
| | flat | 0.50 | 1.14 | 1.00 | 0.96 | 1.18 | 1.00 | 0.99 | 1.17 | 1.00 | 0.98 |
| cycle 1 | steep | 0.25 | 3.98 | 0.99 | 1.67 | 4.80 | 0.97 | 2.02 | 5.14 | 1.02 | 2.14 |
| | flat | 0.50 | 2.39 | 0.97 | 1.29 | 2.47 | 0.98 | 1.34 | 2.74 | 1.01 | 1.48 |
| cycle 4 | steep | 0.25 | 1.04 | 1.00 | 0.97 | 1.12 | 1.03 | 1.05 | 1.19 | 0.99 | 1.12 |
| | flat | 0.50 | 1.04 | 1.00 | 0.97 | 1.04 | 0.99 | 0.98 | 1.18 | 0.92 | 1.11 |
| white | steep | 0.25 | 1.06 | 1.00 | 1.00 | 1.04 | 1.01 | 0.98 | 1.04 | 1.00 | 0.97 |
| noise | flat | 0.50 | 1.06 | 1.00 | 1.00 | 1.05 | 1.00 | 0.99 | 1.02 | 1.00 | 0.96 |

Numbers greater than 1 favor EPSE. Based on logistic model post-stratification in 1000 replications of simple random sampling of size $n = 200$ from a fixed population of size $N = 1000$. Replications in which at least one stratum had fewer than two samples are omitted from the summary: 42 reps for six strata and 2 reps for four strata on steep curve, $\sigma = 0.25$.

Overall, these simulation experiments demonstrate that the practical effect of first fitting a parametric model with a fixed number of parameters to the survey data before post-stratifying is small, even for relatively small sample sizes. Given the types of models actually being used for classification in forest inventory applications, it would be of interest to study the case of semiparametric and nonparametric classification models, either analytically or via simulation, to assess the effect of a large number of unknown parameters that grows with sample size.

## APPENDIX

We begin with design results, in which all probability computations are with respect to the sampling design and the covariate model, specified in D1. To establish the design consistency of the EPSE, we begin with two lemmas.

LEMMA A.1. *Under* D4 *and* D5, *the Horvitz–Thompson estimator*
$$\bar{y}_\pi = \frac{1}{N} \sum_{i \in U_N} \frac{y_i I_{\{i \in s\}}}{\pi_i}$$
*is mean square consistent in the sense that*
$$\mathrm{E}[(\bar{y}_\pi - \bar{y}_N)^2] \to 0 \qquad as\ N \to \infty$$
*and hence design-consistent in the sense that for all $\varepsilon > 0$,*
$$\Pr\{|\bar{y}_\pi - \bar{y}_N| > \varepsilon\} \to 0$$
*as $N \to \infty$.*

PROOF. It suffices to show mean square consistency. Because the Horvitz–Thompson estimator is unbiased,
$$\begin{aligned}
&\mathrm{E}[(\bar{y}_\pi - \bar{y}_N)^2] \\
&= N^{-2} \sum_{i,j \in U_N} \Delta_{ij} \frac{y_i y_j}{\pi_i \pi_j} \\
&\leq \frac{1}{N \pi_N^*} \sum_{i \in U_N} \frac{y_i^2}{N} + \frac{1}{N^2 (\pi_N^*)^2} \left( \sum_{i,j \in U_N : i \neq j} \Delta_{ij}^2 \right)^{1/2} \left( \sum_{i,j \in U_N : i \neq j} y_i^2 y_j^2 \right)^{1/2} \\
&\leq \frac{1}{N \pi_N^*} \sum_{i \in U_N} \frac{y_i^2}{N} \\
&\quad + \frac{1}{N^{1/2+\kappa} (\pi_N^*)^2} \left( \frac{N \max_{i \in U_N} \sum_{j \in U_N : j \neq i} \Delta_{ij}^2}{N^{1-2\kappa}} \right)^{1/2} \left\{ \left( \sum_{i \in U_N} \frac{y_i^2}{N} \right)^2 \right\}^{1/2}
\end{aligned}$$
which converges to zero as $N \to \infty$ under D4 and D5. □

LEMMA A.2. *Assume that* D1–D5 *hold and fix* $h \in \{1, 2, \ldots, H\}$. *Define*

$$L_i(\mathbf{u}) = \begin{cases} 1, & \text{if } m^{-1}(\tau_{h-1}) < \mathbf{x}_i'(\boldsymbol{\lambda} + \mathbf{u}) \leq m^{-1}(\tau_h), \\ 0, & \text{otherwise,} \end{cases}$$

*and* $L_i = L_i(\mathbf{0})$. *Define* $Q_{N1}(\mathbf{u}) = N^{-1} \sum_{i \in U_N} (L_i(\mathbf{u}) - L_i)$, $Q_{N2}(\mathbf{u}) = N^{-1} \times \sum_{i \in U_N} I_{\{i \in s\}} \pi_i^{-1} (L_i(\mathbf{u}) - L_i)$ *and* $Q_{N3}(\mathbf{u}) = N^{-1} \sum_{i \in U_N} y_i I_{\{i \in s\}} \pi_i^{-1} (L_i(\mathbf{u}) - L_i)$. *Then for* $\ell = 1, 2, 3$ *and for all* $\varepsilon > 0$, $\Pr\{|Q_{N\ell}(\hat{\boldsymbol{\lambda}} - \boldsymbol{\lambda})| > \varepsilon\} \to 0$ *as* $N \to \infty$.

PROOF. Consider $Q_{N3}(\mathbf{u})$; the arguments are similar for $Q_{N1}(\mathbf{u})$ and $Q_{N2}(\mathbf{u})$. Let $\varepsilon, \delta > 0$ be given. Then

$$\begin{aligned}
\Pr\{|Q_{N3}(\hat{\boldsymbol{\lambda}} - \boldsymbol{\lambda})| > \varepsilon\} \\
= \Pr\{|Q_{N3}(\hat{\boldsymbol{\lambda}} - \boldsymbol{\lambda})| > \varepsilon, \|\hat{\boldsymbol{\lambda}} - \boldsymbol{\lambda}\| > \delta\} \\
+ \Pr\{|Q_{N3}(\hat{\boldsymbol{\lambda}} - \boldsymbol{\lambda})| > \varepsilon, \|\hat{\boldsymbol{\lambda}} - \boldsymbol{\lambda}\| \leq \delta\} \\
\leq \Pr\{\|\hat{\boldsymbol{\lambda}} - \boldsymbol{\lambda}\| > \delta\} + \Pr\left\{\sup_{\mathbf{u}: \|\mathbf{u}\| \leq \delta} |Q_{N3}(\mathbf{u})| > \varepsilon\right\}.
\end{aligned} \quad (15)$$

By D3, the first term converges to zero as $N \to \infty$. Consider the second term:

$$(16) \quad \Pr\left\{\sup_{\mathbf{u}: \|\mathbf{u}\| \leq \delta} |Q_{N3}(\mathbf{u})| > \varepsilon\right\} \leq \frac{\mathrm{E}[\sup_{\mathbf{u}: \|\mathbf{u}\| \leq \delta} |Q_{N3}(\mathbf{u})|]}{\varepsilon}.$$

Now, using the fact that $|\mathbf{x}_i' \mathbf{u}| \leq M\delta$ from D1,

$$\begin{aligned}
|L_i(\mathbf{u}) - L_i| &= |-I_{\{m^{-1}(\tau_{h-1}) < \mathbf{x}_i' \boldsymbol{\lambda} \leq m^{-1}(\tau_{h-1}) - \mathbf{x}_i' \mathbf{u}\}} + I_{\{m^{-1}(\tau_{h-1}) - \mathbf{x}_i' \mathbf{u} < \mathbf{x}_i' \boldsymbol{\lambda} \leq m^{-1}(\tau_{h-1})\}} \\
&\quad + I_{\{m^{-1}(\tau_h) < \mathbf{x}_i' \boldsymbol{\lambda} \leq m^{-1}(\tau_h) - \mathbf{x}_i' \mathbf{u}\}} - I_{\{m^{-1}(\tau_h) - \mathbf{x}_i' \mathbf{u} < \mathbf{x}_i' \boldsymbol{\lambda} \leq m^{-1}(\tau_h)\}}| \\
&\leq |I_{\{m^{-1}(\tau_{h-1}) < \mathbf{x}_i' \boldsymbol{\lambda} \leq m^{-1}(\tau_{h-1}) + M\delta\}} + I_{\{m^{-1}(\tau_{h-1}) - M\delta < \mathbf{x}_i' \boldsymbol{\lambda} \leq m^{-1}(\tau_{h-1})\}} \\
&\quad + I_{\{m^{-1}(\tau_h) < \mathbf{x}_i' \boldsymbol{\lambda} \leq m^{-1}(\tau_h) + M\delta\}} + I_{\{m^{-1}(\tau_h) - M\delta < \mathbf{x}_i' \boldsymbol{\lambda} \leq m^{-1}(\tau_h)\}}| \\
&=: I_{1i} + I_{2i} + I_{3i} + I_{4i},
\end{aligned}$$

which does not depend on $\mathbf{u}$. Hence

$$\begin{aligned}
\mathrm{E}\left[\sup_{\mathbf{u}: \|\mathbf{u}\| \leq \delta} |Q_{N3}(\mathbf{u})|\right] &\leq \mathrm{E}\left[\sup_{\mathbf{u}: \|\mathbf{u}\| \leq \delta} N^{-1} \sum_{i \in U_N} \frac{|y_i| I_{\{i \in s\}}}{\pi_i} |L_i(\mathbf{u}) - L_i|\right] \\
&\leq \mathrm{E}\left[\sup_{\mathbf{u}: \|\mathbf{u}\| \leq \delta} N^{-1} \sum_{i \in U_N} \frac{|y_i| I_{\{i \in s\}}}{\pi_i} (I_{1i} + I_{2i} + I_{3i} + I_{4i})\right] \\
&= \mathrm{E}\left[N^{-1} \sum_{i \in U_N} \frac{|y_i| I_{\{i \in s\}}}{\pi_i} (I_{1i} + I_{2i} + I_{3i} + I_{4i})\right]
\end{aligned}$$

$$= N^{-1} \sum_{i \in U_N} |y_i|(I_{1i} + I_{2i} + I_{3i} + I_{4i})$$

$$\leq \left(4N^{-1} \sum_{i \in U_N} y_i^2\right)^{1/2} \left(N^{-1} \sum_{j=1}^{4} \sum_{i \in U_N} I_{ji}^2\right)^{1/2}.$$

By D5, it suffices to show that the second term in this product converges to zero as $N \to \infty$. Now

$$N^{-1} \sum_{i \in U_N} I_{1i}^2 = N^{-1} \sum_{i \in U_N} I_{1i}$$

$$= G_{N\boldsymbol{\lambda}}(m^{-1}(\tau_{h-1}) + M\delta) - G_{N\boldsymbol{\lambda}}(m^{-1}(\tau_{h-1}))$$

$$= G_{N\boldsymbol{\lambda}}(m^{-1}(\tau_{h-1}) + M\delta) - G_{\boldsymbol{\lambda}}(m^{-1}(\tau_{h-1}) + M\delta)$$

$$+ G_{\boldsymbol{\lambda}}(m^{-1}(\tau_{h-1}) + M\delta) - G_{\boldsymbol{\lambda}}(m^{-1}(\tau_{h-1}))$$

$$+ G_{\boldsymbol{\lambda}}(m^{-1}(\tau_{h-1})) - G_{N\boldsymbol{\lambda}}(m^{-1}(\tau_{h-1}))$$

$$\leq 2 \sup_{z} |G_{N\boldsymbol{\lambda}}(z) - G_{\boldsymbol{\lambda}}(z)|$$

$$+ \{G_{\boldsymbol{\lambda}}(m^{-1}(\tau_{h-1}) + M\delta) - G_{\boldsymbol{\lambda}}(m^{-1}(\tau_{h-1}))\}.$$

The sup term goes to zero as $N \to \infty$ by D1, and the remaining term in curly braces can be made arbitrarily small because $\delta > 0$ was arbitrary and $m^{-1}(\tau_{h-1})$ is a continuity point of $G_{\boldsymbol{\lambda}}$. Similar arguments can be applied for $I_{2i}$, $I_{3i}$ and $I_{4i}$. It therefore follows that $\Pr\{\sup_{\mathbf{u}:\|\mathbf{u}\| \leq \delta} |Q_{N3}(\mathbf{u})| > \varepsilon\} \to 0$ as $N \to \infty$, so the desired result follows from (15). □

PROOF OF RESULT 1, DESIGN CONSISTENCY OF EPSE. Define
$$\hat{L}_h = \{\mathbf{x}_i : m^{-1}(\tau_{h-1}) < \mathbf{x}_i'\hat{\boldsymbol{\lambda}} \leq m^{-1}(\tau_h)\}$$
and define $L_h$ similarly, with $\hat{\boldsymbol{\lambda}}$ replaced by $\boldsymbol{\lambda}$. For fixed $h \in \{1, 2, \ldots, H\}$, define

$$F_h = G_{\boldsymbol{\lambda}}(m^{-1}(\tau_h)) - G_{\boldsymbol{\lambda}}(m^{-1}(\tau_{h-1})), \qquad g(w_1, w_2) = \frac{w_1}{w_2 + F_h},$$

$$\hat{W}_{1Nh} = N^{-1} \sum_{i \in U_N} I_{\{\mathbf{x}_i \in \hat{L}_h\}}, \qquad W_{1Nh} = N^{-1} \sum_{i \in U_N} I_{\{\mathbf{x}_i \in L_h\}},$$

$$\hat{W}_{2Nh} = N^{-1} \sum_{i \in U_N} I_{\{\mathbf{x}_i \in \hat{L}_h\}} I_{\{i \in s\}} \pi_i^{-1} - F_h, \qquad W_{2Nh} = N^{-1} \sum_{i \in U_N} I_{\{\mathbf{x}_i \in L_h\}} - F_h,$$

$$\hat{W}_{3Nh} = N^{-1} \sum_{i \in U_N} y_i I_{\{\mathbf{x}_i \in \hat{L}_h\}} I_{\{i \in s\}} \pi_i^{-1}, \qquad W_{3Nh} = N^{-1} \sum_{i \in U_N} y_i I_{\{\mathbf{x}_i \in L_h\}}.$$

Note that $g(\cdot)$ is continuous for $w_2 \neq -F_h$ and that $\hat{W}_{1Nh}$ and $W_{1Nh}$ are bounded by 1. Choose $\delta \in (0, F_h)$, where $F_h > 0$ by D2. Then

$$\Pr\{|W_{2Nh}| > \delta\} = \Pr\{|G_{N\boldsymbol{\lambda}}(m^{-1}(\tau_h)) - G_{N\boldsymbol{\lambda}}(m^{-1}(\tau_{h-1})) - F_h| > \delta\},$$

which goes to zero as $N \to \infty$ by D1. Combining Lemma A.1 and Lemma A.2 we have that for all $\delta > 0$, $\Pr\{|\hat{W}_{kNh} - W_{kNh}| > \delta\} \to 0$ as $N \to \infty$ for $k = 1, 2, 3$, so

$$\Pr\{|\hat{W}_{2Nh}| > \delta\} \leq \Pr\{|\hat{W}_{2Nh} - W_{2Nh}| > \delta/2\} + \Pr\{|W_{2Nh}| > \delta/2\} \to 0$$

as $N \to \infty$. Now given any $\varepsilon > 0$,

$$\Pr\{|g(\hat{W}_{1Nh}, \hat{W}_{2Nh})\hat{W}_{3Nh} - g(W_{1Nh}, W_{2Nh})W_{3Nh}| > \varepsilon\}$$
$$\leq \Pr\{|\hat{W}_{2Nh}| > \delta\} + \Pr\{|W_{2Nh}| > \delta\}$$
$$+ \Pr\{|g(\hat{W}_{1Nh}, \hat{W}_{2Nh})\hat{W}_{3Nh} - g(W_{1Nh}, W_{2Nh})W_{3Nh}| > \varepsilon,$$
$$|\hat{W}_{2Nh}| \leq \delta, |W_{2Nh}| \leq \delta\}.$$

From above, the first and second probabilities go to zero as $N \to \infty$, so consider the third term. Write $\hat{g} = g(\hat{W}_{1Nh}, \hat{W}_{2Nh})$ and $g = g(W_{1Nh}, W_{2Nh})$ and note that

$$\hat{g}\hat{W}_{3Nh} - gW_{3Nh} = \hat{g}(\hat{W}_{3Nh} - W_{3Nh}) + (\hat{g} - g)W_{3Nh},$$

so that the third term in the probability statement above is bounded by

$$\Pr\{|\hat{g}||\hat{W}_{3Nh} - W_{3Nh}| > \varepsilon/2, |\hat{W}_{2Nh}| \leq \delta\}$$
$$+ \Pr\{|\hat{g} - g||W_{3Nh}| > \varepsilon/2, |\hat{W}_{2Nh}| \leq \delta, |W_{2Nh}| \leq \delta\}.$$

Now $g(w_1, w_2)$ is uniformly continuous and bounded between zero and $(F_h - \delta)^{-1}$ on the set where $|w_2| \leq \delta < F_h$, so the first term converges to zero as $N \to \infty$ by Lemma A.2. For the second term, first note that $|W_{3Nh}| \leq N^{-1}\sum_{i \in U_N} |y_i| \leq (N^{-1}\sum_{i \in U_N} y_i^2)^{1/2}$, which is finite by D5. Then, by uniform continuity of $g$ there exists $\delta_\varepsilon > 0$ such that for all $N$,

$$\Pr\{|g(\hat{W}_{1Nh}, \hat{W}_{2Nh}) - g(W_{1Nh}, W_{2Nh})||W_{3Nh}| > \varepsilon/2, |\hat{W}_{2Nh}| \leq \delta, |W_{2Nh}| \leq \delta\}$$
$$\leq \Pr\{\|(\hat{W}_{1Nh}, \hat{W}_{2Nh}) - (W_{1Nh}, W_{2Nh})\| > \delta_\varepsilon\},$$

which converges to zero as $N \to \infty$ by Lemma A.2. Since the above results hold for each $h \in \{1, 2, \ldots, H\}$, it follows that

$$\Pr\left\{\left|\hat{\mu}_y^*(\hat{\boldsymbol{\lambda}}) - N^{-1}\sum_{i \in U_N} y_i\right| > \varepsilon\right\}$$

(17)
$$= \Pr\left\{\left|\sum_{h=1}^{H}\left(\frac{\hat{W}_{1Nh}\hat{W}_{3Nh}}{\hat{W}_{2Nh} + F_h} - \frac{W_{1Nh}W_{3Nh}}{W_{2Nh} + F_h}\right)\right| > \varepsilon\right\}$$
$$\to 0 \quad \text{as } N \to \infty,$$

and the result is proved. $\square$

PROOF OF RESULT 2, CENTRAL LIMIT THEOREM FOR EPSE UNDER SUPERPOPULATION MODEL. In what follows, the probability mechanism is the joint distribution for $(\mathbf{x}_i', y_i)$ as given in model assumptions M1, M2 and M4. In particular, all expectation, probability and order in probability statements are with respect to this superpopulation model. Let $K(\boldsymbol{\lambda})$ be a neighborhood of $\boldsymbol{\lambda}$ which is bounded away from $\mathbf{0}$, and consider $\boldsymbol{\gamma} \in K(\boldsymbol{\lambda})$. For $\ell \in \{0, 1, 2\}$, both $A_{Nh\ell}(\boldsymbol{\gamma})$ and $A_{nh\ell}(\boldsymbol{\gamma})$ are $U$-statistics with kernel $y_i^\ell I_{\{\tau_{h-1} < m(\boldsymbol{\gamma}'\mathbf{x}_i) \le \tau_h\}}$ and common expectation $\alpha_{h\ell}(\boldsymbol{\gamma})$ given in (9). We will apply Theorem 2.8 of Randles [8] to derive asymptotic approximations for $A_{Nh\ell}(\boldsymbol{\gamma})$ and $A_{nh\ell}(\boldsymbol{\gamma})$ for $\boldsymbol{\gamma} \in K(\boldsymbol{\lambda})$, which requires checking that Conditions 2.2 and 2.3 on page 465 of Randles [8] hold.

Condition 2.2 is immediate by M3. To verify Condition 2.3, we need to establish that (2.4) and (2.5) on page 465 of Randles [8] hold. Let $D(\boldsymbol{\gamma}, d) \subset K(\boldsymbol{\lambda})$ be a sphere of radius $d$ centered at $\boldsymbol{\gamma}$, and let $\boldsymbol{\theta} \in D(\boldsymbol{\gamma}, d)$. Then the supremum random variable in (2.4) of Randles [8] is

$$
\begin{aligned}
&\sup_{\boldsymbol{\theta} \in D(\boldsymbol{\gamma},d)} |y_i^\ell I_{\{\tau_{h-1}<m(\boldsymbol{\theta}'\mathbf{x}_i) \le \tau_h\}} - y_i^\ell I_{\{\tau_{h-1}<m(\boldsymbol{\gamma}'\mathbf{x}_i) \le \tau_h\}}| \\
&= |y_i|^\ell \sup_{\boldsymbol{\theta} \in D(\boldsymbol{\gamma},d)} |I_{\{\tau_{h-1}<m(\boldsymbol{\theta}'\mathbf{x}_i) \le \tau_h\}} - I_{\{\tau_{h-1}<m(\boldsymbol{\gamma}'\mathbf{x}_i) \le \tau_h\}}|.
\end{aligned}
$$
(18)

Since $\mathrm{E}(|y_i|^\ell | \mathbf{x}_i) \le 1 + K_1^{1/4} + K_1^{1/2} + K_1 =: K_2$ for $\ell = 0, 1, 2$ and 4 by hypothesis, it suffices to look at the difference of the indicators in (18). Similarly to the reasoning in the proof of Lemma A.2,

$$
\begin{aligned}
&|I_{\{\tau_{h-1}<m(\boldsymbol{\theta}'\mathbf{x}_i) \le \tau_h\}} - I_{\{\tau_{h-1}<m(\boldsymbol{\gamma}'\mathbf{x}_i) \le \tau_h\}}| \\
&\le I_{\{\boldsymbol{\theta}'\mathbf{x}_i < m^{-1}(\tau_h) < \boldsymbol{\gamma}'\mathbf{x}_i\}} + I_{\{\boldsymbol{\gamma}'\mathbf{x}_i < m^{-1}(\tau_{h-1}) < \boldsymbol{\theta}'\mathbf{x}_i\}} \\
&\quad + I_{\{\boldsymbol{\gamma}'\mathbf{x}_i < m^{-1}(\tau_h) < \boldsymbol{\theta}'\mathbf{x}_i\}} + I_{\{\boldsymbol{\theta}'\mathbf{x}_i < m^{-1}(\tau_{h-1}) < \boldsymbol{\gamma}'\mathbf{x}_i\}}.
\end{aligned}
$$
(19)

We now bound this sum by maximizing the indicated events with respect to $\boldsymbol{\theta}$, subject to the constraint that $\boldsymbol{\theta}$ is in the closure of $D(\boldsymbol{\gamma}, d)$. Let $\boldsymbol{\theta}_{\inf} = \arg\inf_{\boldsymbol{\theta} \in D(\boldsymbol{\gamma},d)} \boldsymbol{\theta}'\mathbf{x}_i$ and $\boldsymbol{\theta}_{\sup} = \arg\sup_{\boldsymbol{\theta} \in D(\boldsymbol{\gamma},d)} \boldsymbol{\theta}'\mathbf{x}_i$. Then, $\boldsymbol{\theta}_{\inf}$ and $\boldsymbol{\theta}_{\sup}$ must occur on the boundary of the sphere $D(\boldsymbol{\gamma}, d)$, since any point on the interior of the sphere has nonzero derivative for the linear function $\boldsymbol{\theta}'\mathbf{x}_i$. Optimizing $\boldsymbol{\theta}'\mathbf{x}_i$ subject to $d^2 = (\boldsymbol{\gamma} - \boldsymbol{\theta})'(\boldsymbol{\gamma} - \boldsymbol{\theta})$, we have that $(\boldsymbol{\theta}_{\inf}, \boldsymbol{\theta}_{\sup}) = (\boldsymbol{\gamma} - d\mathbf{x}_i/\|\mathbf{x}_i\|, \boldsymbol{\gamma} + d\mathbf{x}_i/\|\mathbf{x}_i\|)$, so that the sum of the indicators in (19) is bounded above by

$$
\begin{aligned}
&I_{\{\boldsymbol{\gamma}'\mathbf{x}_i - d\|\mathbf{x}_i\| \le m^{-1}(\tau_h) \le \boldsymbol{\gamma}'\mathbf{x}_i\}} + I_{\{\boldsymbol{\gamma}'\mathbf{x}_i \le m^{-1}(\tau_{h-1}) \le \boldsymbol{\gamma}'\mathbf{x}_i + d\|\mathbf{x}_i\|\}} \\
&\quad + I_{\{\boldsymbol{\gamma}'\mathbf{x}_i \le m^{-1}(\tau_h) \le \boldsymbol{\gamma}'\mathbf{x}_i + d\|\mathbf{x}_i\|\}} + I_{\{\boldsymbol{\gamma}'\mathbf{x}_i - d\|\mathbf{x}_i\| \le m^{-1}(\tau_{h-1}) \le \boldsymbol{\gamma}'\mathbf{x}_i\}}.
\end{aligned}
$$
(20)

Consider the first of these four indicators. Define

$$G_{\boldsymbol{\gamma}}(t) = \Pr\{\boldsymbol{\gamma}'\mathbf{x}_i \le t\}$$

$$= \int_{-\infty}^{t} \frac{1}{|\gamma_1|} \bigg\{ \int \cdots \int f((s - \gamma_2 x_2 - \cdots - \gamma_p x_p)/\gamma,$$
$$x_2, \ldots, x_p) \, dx_2 \cdots dx_p \bigg\} \, ds$$

(assume without loss of generality that $\gamma_1^*$, the first element of $\boldsymbol{\gamma}$, is not zero). Since $\|\mathbf{x}_i\| \leq M < \infty$ with probability 1,

$$\mathrm{E}(I_{\{\boldsymbol{\gamma}'\mathbf{x}_i - d\|\mathbf{x}_i\| \leq m^{-1}(\tau_h) \leq \boldsymbol{\gamma}'\mathbf{x}_i\}}) \leq \mathrm{E}(I_{\{\boldsymbol{\gamma}'\mathbf{x}_i - dM \leq m^{-1}(\tau_h) \leq \boldsymbol{\gamma}'\mathbf{x}_i\}})$$
$$= G_{\boldsymbol{\gamma}}(m^{-1}(\tau_h) + dM) - G_{\boldsymbol{\gamma}}(m^{-1}(\tau_h))$$
$$\leq dK_3$$

for some constant $K_3$, using the mean value theorem and the continuity and compact support of $f$. Arguing in this fashion for the remaining three indicators in (20) establishes (2.4) of Randles [8].

Next,

$$\lim_{d \to 0} \mathrm{E}\bigg[\sup_{\boldsymbol{\theta} \in D(\boldsymbol{\gamma},d)} |y_i^\ell I_{\{\tau_{h-1} < m(\boldsymbol{\theta}'\mathbf{x}_i) \leq \tau_h\}} - y_i^\ell I_{\{\tau_{h-1} < m(\boldsymbol{\gamma}'\mathbf{x}_i) \leq \tau_h\}}|^2\bigg]$$

$$= \lim_{d \to 0} \mathrm{E}\bigg[y_i^{2\ell} \sup_{\boldsymbol{\theta} \in D(\boldsymbol{\gamma},d)} |I_{\{\tau_{h-1} < m(\boldsymbol{\theta}'\mathbf{x}_i) \leq \tau_h\}} - I_{\{\tau_{h-1} < m(\boldsymbol{\gamma}'\mathbf{x}_i) \leq \tau_h\}}|\bigg]$$

$$\leq \lim_{d \to 0} K_2 \mathrm{E}[I_{\{\boldsymbol{\gamma}'\mathbf{x}_i - d\|\mathbf{x}_i\| \leq m^{-1}(\tau_h) \leq \boldsymbol{\gamma}'\mathbf{x}_i\}} + I_{\{\boldsymbol{\gamma}'\mathbf{x}_i \leq m^{-1}(\tau_{h-1}) \leq \boldsymbol{\gamma}'\mathbf{x}_i + d\|\mathbf{x}_i\|\}}$$
$$+ I_{\{\boldsymbol{\gamma}'\mathbf{x}_i \leq m^{-1}(\tau_h) \leq \boldsymbol{\gamma}'\mathbf{x}_i + d\|\mathbf{x}_i\|\}} + I_{\{\boldsymbol{\gamma}'\mathbf{x}_i - d\|\mathbf{x}_i\| \leq m^{-1}(\tau_{h-1}) \leq \boldsymbol{\gamma}'\mathbf{x}_i\}}]$$

$$= 0$$

by the previous linear bound on the expectation, so that (2.5) of Randles [8] is satisfied. It follows from Randles' Theorem 2.8 that

(21) $\qquad A_{Nh\ell}(\hat{\boldsymbol{\lambda}}) = \alpha_{h\ell}(\hat{\boldsymbol{\lambda}}) + A_{Nh\ell}(\boldsymbol{\lambda}) - \alpha_{h\ell}(\boldsymbol{\lambda}) + o_p(N^{-1/2}),$

(22) $\qquad A_{nh\ell}(\hat{\boldsymbol{\lambda}}) = \alpha_{h\ell}(\hat{\boldsymbol{\lambda}}) + A_{nh\ell}(\boldsymbol{\lambda}) - \alpha_{h\ell}(\boldsymbol{\lambda}) + o_p(n^{-1/2}).$

Define $a_h = A_{Nh0}(\boldsymbol{\lambda}) - A_{nh0}(\boldsymbol{\lambda})$ and $b_h = A_{Nh1}(\boldsymbol{\lambda}) - A_{nh1}(\boldsymbol{\lambda})$. Straightforward calculations show that

$$\mathrm{Cov}(a_h, a_k) = \frac{1}{n}\bigg(1 - \frac{n}{N}\bigg)(\alpha_{h0}(\boldsymbol{\lambda})I_{\{h=k\}} - \alpha_{h0}(\boldsymbol{\lambda})\alpha_{k0}(\boldsymbol{\lambda})),$$

$$\mathrm{Cov}(a_h, b_k) = \frac{1}{n}\bigg(1 - \frac{n}{N}\bigg)(\alpha_{h1}(\boldsymbol{\lambda})I_{\{h=k\}} - \alpha_{h0}(\boldsymbol{\lambda})\alpha_{k1}(\boldsymbol{\lambda})),$$

from which it follows that $a_h = O_p(n^{-1/2})$ and $b_h = O_p(n^{-1/2})$. Also note that $\alpha_{h\ell}(\hat{\boldsymbol{\lambda}}) - \alpha_{h\ell}(\boldsymbol{\lambda}) = o_p(1)$ by M3 and M4, and that

$$A_{Nh\ell}(\boldsymbol{\lambda}) - \alpha_{h\ell}(\boldsymbol{\lambda}) = O_p(N^{-1/2}) \quad \text{and} \quad A_{nh\ell}(\boldsymbol{\lambda}) - \alpha_{h\ell}(\boldsymbol{\lambda}) = O_p(n^{-1/2})$$

by the central limit theorem.

Since $\bar{y}_N = \sum_{h=1}^{H} A_{Nh1}(\gamma)$ for any $\gamma$, we have that

$$(23) \qquad \hat{\mu}_y(\hat{\boldsymbol{\lambda}}) - \bar{y}_N = \sum_{h=1}^{H} \left\{ \frac{A_{Nh0}(\hat{\boldsymbol{\lambda}})A_{nh1}(\hat{\boldsymbol{\lambda}}) - A_{nh0}(\hat{\boldsymbol{\lambda}})A_{Nh1}(\hat{\boldsymbol{\lambda}})}{A_{nh0}(\hat{\boldsymbol{\lambda}})} \right\}.$$

Substituting (21) and (22) in the numerator of the summand in (23), we apply the order results above to obtain

$$\begin{aligned}
A_{Nh0}(\hat{\boldsymbol{\lambda}})&A_{nh1}(\hat{\boldsymbol{\lambda}}) - A_{nh0}(\hat{\boldsymbol{\lambda}})A_{Nh1}(\hat{\boldsymbol{\lambda}}) \\
&= (\alpha_{h0}(\hat{\boldsymbol{\lambda}}) - \alpha_{h0}(\boldsymbol{\lambda}))(A_{nh1}(\boldsymbol{\lambda}) - A_{Nh1}(\boldsymbol{\lambda})) \\
&\quad - (\alpha_{h1}(\hat{\boldsymbol{\lambda}}) - \alpha_{h1}(\boldsymbol{\lambda}))(A_{nh0}(\boldsymbol{\lambda}) - A_{Nh0}(\boldsymbol{\lambda})) \\
&\quad + A_{Nh0}(\boldsymbol{\lambda})A_{nh1}(\boldsymbol{\lambda}) - A_{nh0}(\boldsymbol{\lambda})A_{Nh1}(\boldsymbol{\lambda}) + o_p(n^{-1/2}) \\
&= A_{Nh0}(\boldsymbol{\lambda})A_{nh1}(\boldsymbol{\lambda}) - A_{nh0}(\boldsymbol{\lambda})A_{Nh1}(\boldsymbol{\lambda}) + o_p(n^{-1/2}) \\
(24) \\
&= (A_{Nh0}(\boldsymbol{\lambda}) - \alpha_{h0}(\boldsymbol{\lambda}) + \alpha_{h0}(\boldsymbol{\lambda}))(A_{nh1}(\boldsymbol{\lambda}) - \alpha_{h1}(\boldsymbol{\lambda}) + \alpha_{h1}(\boldsymbol{\lambda})) \\
&\quad - (A_{nh0}(\boldsymbol{\lambda}) - \alpha_{h0}(\boldsymbol{\lambda}) + \alpha_{h0}(\boldsymbol{\lambda})) \\
&\quad\quad \times (A_{Nh1}(\boldsymbol{\lambda}) - \alpha_{h1}(\boldsymbol{\lambda}) + \alpha_{h1}(\boldsymbol{\lambda})) + o_p(n^{-1/2}) \\
&= \alpha_{h1}(\boldsymbol{\lambda})(A_{Nh0}(\boldsymbol{\lambda}) - A_{nh0}(\boldsymbol{\lambda})) \\
&\quad + \alpha_{h0}(\boldsymbol{\lambda})(A_{nh1}(\boldsymbol{\lambda}) - A_{Nh1}(\boldsymbol{\lambda})) + o_p(n^{-1/2}),
\end{aligned}$$

where we have used the facts that $A_{nhl}(\boldsymbol{\lambda})$ and $A_{Nhl}(\boldsymbol{\lambda})$ are $O_p(1)$ by the weak law of large numbers.

From (22), the denominator of the summand in (23) is $A_{nh0}(\hat{\boldsymbol{\lambda}}) = \alpha_{h0}(\boldsymbol{\lambda}) + o_p(1)$, and so

$$(25) \qquad \frac{1}{A_{nh0}(\hat{\boldsymbol{\lambda}})} = \frac{1}{\alpha_{h0}(\boldsymbol{\lambda})} + o_p(1)$$

since $\alpha_{h0}(\boldsymbol{\lambda}) > 0$ by M4.

Substituting (24) and (25) into (23), we have

$$\hat{\mu}_y(\hat{\boldsymbol{\lambda}}) - \bar{y}_N = \sum_{h=1}^{H} \left\{ \frac{\alpha_{h1}(\boldsymbol{\lambda})}{\alpha_{h0}(\boldsymbol{\lambda})} (A_{Nh0}(\boldsymbol{\lambda}) - A_{nh0}(\boldsymbol{\lambda})) - (A_{Nh1}(\boldsymbol{\lambda}) - A_{nh1}(\boldsymbol{\lambda})) \right\}$$
$$(26)$$
$$\qquad + o_p(n^{-1/2}),$$

so that the asymptotic distribution is the same as that obtained when $\boldsymbol{\lambda}$ is known.

It remains to derive this asymptotic distribution. Applying the central limit theorem to (26), we have that the limiting distribution of the EPSE error is normal with mean zero. Using earlier covariance computations, and the fact that $\sum_{h=1}^{H} b_h = \bar{y}_N - \bar{y}_\pi$, it follows that the variance of the leading terms in (26) is approximated by

$$
\begin{aligned}
(27) \quad & \operatorname{Var}(\hat{\mu}_y(\hat{\boldsymbol{\lambda}}) - \bar{y}_N) \\
& \simeq -\frac{1}{n}\left(1 - \frac{n}{N}\right) \sum_{h=1}^{H} \frac{\alpha_{h1}^2(\boldsymbol{\lambda})}{\alpha_{h0}(\boldsymbol{\lambda})} + \frac{1}{n}\left(1 - \frac{n}{N}\right)\left(\sum_{h=1}^{H} \alpha_{h1}(\boldsymbol{\lambda})\right)^2 \\
& \quad + \operatorname{Var}(\bar{y}_\pi - \bar{y}_N) \\
& = \frac{1}{n}\left(1 - \frac{n}{N}\right)\left\{-\sum_{h=1}^{H} \frac{\alpha_{h1}^2(\boldsymbol{\lambda})}{\alpha_{h0}(\boldsymbol{\lambda})} + [\mathrm{E}(y_i)]^2 + \operatorname{Var}(y_i)\right\}.
\end{aligned}
$$

Note that, by definition of expectation given an event,

$$\frac{\alpha_{h1}(\boldsymbol{\lambda})}{\alpha_{h0}(\boldsymbol{\lambda})} = \mathrm{E}(y_i | \tau_{h-1} < m(\boldsymbol{\lambda}'\mathbf{x}_i) \leq \tau_h)$$

and

$$
\begin{aligned}
\mathrm{E}(y_i^2) = \sum_{h=1}^{H} \alpha_{h0}(\boldsymbol{\lambda})\{&\operatorname{Var}(y_i | \tau_{h-1} < m(\boldsymbol{\lambda}'\mathbf{x}_i) \leq \tau_h) \\
& + [\mathrm{E}(y_i | \tau_{h-1} < m(\boldsymbol{\lambda}'\mathbf{x}_i) \leq \tau_h)]^2\}
\end{aligned}
$$

from which the variance given in Result 2 immediately follows. $\square$

PROOF OF RESULT 3, CONSISTENT ESTIMATION OF EPSE VARIANCE UNDER SUPERPOPULATION MODEL. Note that $A_{Nh\ell}(\boldsymbol{\lambda}) \xrightarrow{p} \alpha_{h\ell}(\boldsymbol{\lambda})$ and $A_{nh\ell}(\boldsymbol{\lambda}) \xrightarrow{p} \alpha_{h\ell}(\boldsymbol{\lambda})$ as $n, N \to \infty$ by the weak law of large numbers, and $\alpha_{h\ell}(\hat{\boldsymbol{\lambda}}) \xrightarrow{p} \alpha_{h\ell}(\boldsymbol{\lambda})$ by continuity of $\alpha_{h\ell}(\cdot)$ for $\ell = 0, 1, 2$. Using (21) and (22) of the Appendix, the term in curly braces in (14) converges in probability to

$$\sum_{h=1}^{H} \alpha_{h0}(\boldsymbol{\lambda})\left\{\frac{\alpha_{h2}(\boldsymbol{\lambda})}{\alpha_{h0}(\boldsymbol{\lambda})} - \left(\frac{\alpha_{h1}(\boldsymbol{\lambda})}{\alpha_{h0}(\boldsymbol{\lambda})}\right)^2\right\}$$

from which the result follows by Slutsky's theorem and Result 2. $\square$

**Acknowledgment.** The authors are grateful to the Associate Editor for suggesting directions in which an earlier version of D4 could be weakened.

DEPARTMENT OF STATISTICS
COLORADO STATE UNIVERSITY
102 STATISTICS BUILDING
FORT COLLINS, COLORADO 80523
USA
E-MAIL: jbreidt@stat.colostate.edu

DEPARTMENT OF STATISTICS
COLORADO STATE UNIVERSITY
201 STATISTICS BUILDING
FORT COLLINS, COLORADO 80523
USA
E-MAIL: jopsomer@stat.colostate.edu